\documentclass[oneside]{article}
\usepackage{amsmath}
\usepackage{amssymb}
\usepackage{bbm}
\usepackage[margin=3cm]{geometry}
\usepackage{hyperref}
\usepackage{tikz}
\usepackage{braket}

\newtheorem{thm}{Theorem}[section]

\newtheorem{lemma}[thm]{Lemma}

\newtheorem{theorem}[thm]{Theorem}
\newtheorem{remark}[thm]{Remark}
\newtheorem{proposition}[thm]{Proposition}

\newtheorem{corollary}[thm]{Corollary}

\newcommand{\R}{\mathbb{R}}
\newcommand{\N}{\mathbb{N}}
\newcommand{\C}{\mathbb{C}}
\newcommand{\E}{\mathbb{E}}
\newcommand{\Z}{\mathbb{Z}}

\renewcommand{\P}{\mathbb{P}}

\title{Two Examples of COM Bounds using Spectral Gaps: Length of the LIS in a Random Permutation and Lipschitz Functions of 1d Markov Chains}
\author{Michael A.~Fr\"olich${}^{1}$ and 
Shannon Starr
${}^2$
\\
\small
${}^{1}$ Department of Anesthesiology and Perioperative Medicine,\\ 
\small School of Medicine, University of
 Alabama at Birmingham (UAB)\\[2pt]
\small ${}^{2}$ 
Department of Applied Mathematics, UAB, Birmingham, AL 35294--1170
}

\date{January 4, 2018}

\begin{document}

\maketitle

\begin{abstract}
\setcounter{section}{0}
We consider two examples for a well-known method for obtaining concentration of measure (COM) bounds for a given observable in a given measure.
The method is to consider an auxiliary Markov chain for which the 
invariant distribution is the measure of interest.
Then one obtains COM bounds involving two quantities.
The first is  
the spectral gap of the Markov transition matrix.
The second is  an appropriate Lipschitz constant for the observable of interest 
with respect to 1 step of the Markov chain.

We consider two examples of the basic method.
The first is to obtain rough COM bounds for the length of the longest increasing subsequence (LIS)
in a uniform random permutation.
The bounds are similar to well-known bounds of Talagrand using his isoperimetric inequality.

The second example is to consider a 1d Markov chain: $X_0,X_1,\dots,X_n$.
We assume the invariant measure for the chain $\mu$ is reversible, and let the initial distribution of $X_0$ be $\mu$.
Then the observable of interest is any function $f(X_0,X_1,\dots,X_n)$, which is Lipschitz with respect to replacement
of single variables.
One  case of this is ``target frequency analysis,''
which is of interest in biostatistics. The auxiliary Markov chain is Glauber dynamics which is gapped in 1d.
\end{abstract}

\section{Statement of the General Method} 

This article is about obtaining concentration of measure (COM) bounds for certain observables in given measures.

For the present article, suppose that $\mathcal{X}$ is a fixed, finite set.
Suppose we are interested in COM bounds for an observable, by which we mean a real-valued function $f$ on $\mathcal{X}$.
And we are interested in COM bounds for the observable $f$ relative to a non-degenerate measure $\mu$ on $\mathcal{X}$.
Let $\mathcal{M}_{+,1}(\mathcal{X})$ denote the set of all probability measures $\mu$ on $\mathcal{X}$.
(So we have $\forall x \in \mathcal{X}$ that $\mu(x)\geq 0$ and we also have $\sum_{x \in \mathcal{X}} \mu(x) = 1$.)
Let $\mathcal{C}(\mathcal{X})$ denote the set of all real-valued functions
$f : \mathcal{X} \to \R$.

Given a pair $(\mu,f)$ from $\mathcal{M}_{+,1}(\mathcal{X})\times \mathcal{C}(\mathcal{X})$, let us define the positive and negative fluctuations as
\begin{equation}
    \mathcal{F}^+_{\mu}(f;a)\, \stackrel{\mathrm{def}}{:=}\, \mu\left(\{x \in \mathcal{X}\, :\, f(x) - \E_{\mu}[f] \geq a\}\right)\, \text{ and }\,
    \mathcal{F}^-_{\mu}(f;a)\, \stackrel{\mathrm{def}}{:=}\, \mathcal{F}^+_{\mu}(-f;a)\, ,
\end{equation}
for $a\in [0,\infty)$.
We are interested in concentration of measure bounds which are bounds on $\mathcal{F}^{\pm}_{\mu}(f;\cdot)$ for a particular pair $(\mu,f)$.
In this article, we will sometimes derive different types of bounds for $\mathcal{F}^+_{\mu}(f;\cdot)$ than for $\mathcal{F}^-_{\mu}(f;\cdot)$.
But if we have sufficiently good bounds for one, then that leads to (possibly weaker) bounds for the other one by Markov's inequality, using the 
fact that $\int_{0}^{\infty} \mathcal{F}^+_{\mu}(f;t)\, dt = \int_{0}^{\infty} \mathcal{F}^-_{\mu}(f;t)\, dt$.

In this setting, Aida and Stroock described a useful elementary method for 
obtaining such bounds, using Chebyshev's inequality \cite{AS}.
They did this along the way to proving even more sophisticated bounds, but for the present article we focus on their first, elementary method.

Suppose one has a Markov chain given by transition matrix $P:\mathcal{X}\times \mathcal{X} \to \R$, having the property that $\mu$ is stationary  for $P$.
So
\begin{equation}
\label{eq:stationarity1}
    \mu \cdot P\, =\, \mu\, ,\ \text{ that is: }\quad \forall x \in \mathcal{X}\, ,\ \text{ we have }\, \mu(x)\, =\, \sum_{y \in \mathcal{X}} \mu(y) P(y,x)\, .
\end{equation}
Also, suppose that moreover $\mu$ is reversible for $P$, meaning
\begin{equation}
\label{eq:reversibility1}
    \forall x,y \in \mathcal{X}\, , \text{ we have }\, \mu(x) P(x,y)\, =\, \mu(y) P(y,x)\, .
\end{equation}
(The condition (\ref{eq:reversibility1}) implies (\ref{eq:stationarity1}), of course,
since $P$ satisfies $P(x,y)\geq 0$ for all $x,y \in \mathcal{X}$ and $\sum_{y\in\mathcal{X}}P(x,y)=1$ for all $x \in \mathcal{X}$.)
Then Aida and Stroock used the Markov chain to derive COM bounds in the measure $\mu$.

Denote the variance in the measure $\mu$ as $\operatorname{Var}_{\mu}(f) = \sum_{x \in \Omega} \big(f(x) - \E_{\mu}[f]\big)^2 \mu(x)$, as usual.
The Dirichlet form for $P$, relative to $\mu$, is denoted by $\mathcal{E}_{P,\mu} : \mathcal{C}(\mathcal{X}) \times \mathcal{C}(\mathcal{X}) \to \R$,
defined as
\begin{equation}
    \mathcal{E}_{P,\mu}(f,g)\, =\, \frac{1}{2}\, \sum_{x,y \in \mathcal{X}} \big(f(x)-f(y)\big)\big(g(x)-g(y)\big)\, \mu(x) P(x,y)\, .
\end{equation}
It is well-known (see for example, Section 13.2 of \cite{LevinPeres}) that
\begin{equation}
\label{eq:DirichletFormula}
    \mathcal{E}_{P,\mu}(f,g)\, =\, \sum_{x \in \mathcal{X}} f(x) \cdot \big((I-P)g\big)(x) \mu(x)\, .
\end{equation}
Since $\mu$ is a reversible measure relative to the Markov transition matrix $P$, the Dirichlet form  $\mathcal{E}_{P,\mu}$ is a positive semi-definite bilinear form on the vector space $\mathcal{C}(\mathcal{X})$.
Then the spectral
gap of $P$, relative to the reversible measure $\mu$, is
\begin{equation}
\label{eq:GapDef}
    \Lambda^{(1)}_{\mu}(P)\, =\, \min\left(\left\{ \mathcal{E}_{P,\mu}(f,f)\, :\, f \in \mathcal{C}(\mathcal{X})\, ,\ \operatorname{Var}_{\mu}(f)\geq 1 \right\}\right)\, 
    =\, \frac{1}{\sup\left(\left\{\operatorname{Var}_{\mu}(f)\, :\, \mathcal{E}_{P,\mu}(f,f)\, \leq\, 1\right\}\right)}\, ,
\end{equation}
which is well-defined as long as $\mu$ is non-degenerate ($\forall x \in \mathcal{X}$, we have $\mu(x)>0$) and $|\mathcal{X}|\neq 1$.
The spectral gap $\Lambda^{(1)}_{\mu}(P)$
is strictly positive as long as $P$ is irreducible ($\forall x,y \in \mathcal{X}$, there is a $T \in \{1,2,\dots\}$
and $x_1,\dots,x_{T-1} \in \mathcal{X}$ such that $P(x_t,x_{t+1})>0$ for all $t \in \{0,\dots,T-1\}$ where $x_0=x$ and $x_T=y$).
See, for example, Levin and Peres for the notation (used in this article) related to finite state space Markov chains \cite{LevinPeres}.

Then, also, let us denote  the $L^{\infty}$-$L^2$ Lipschitz constant of $f \in \mathcal{C}(\mathcal{X})$ with respect to 1 step of the Markov chain
\begin{equation}
\label{eq:LipDef}
    \Phi^{\mathrm{Lip}}_P(f)\, \stackrel{\mathrm{def}}{:=}\, \max_{x \in \mathcal{X}} \left(\sum_{y \in \mathcal{X}} P(x,y) |f(x)-f(y)|^2\right)^{1/2}\, .
\end{equation}
The Aida-Stroock theorem then gives exponential moment generating function bounds.
\begin{theorem}
\label{thm:AS}
Define a function $\Theta : [0,1) \to [0,\infty)$ as
\begin{equation}
    \Theta(t)\, =\, \sum_{n=0}^{\infty} 2^n\, \ln\left(\frac{1}{1-t\cdot 2^{-2n}}\right)\, .
\end{equation}
Then, for $\lambda \geq 0$, satisfying $\lambda < 2 \left(\Lambda^{(1)}_P(\mu)\right)^{1/2}/\Phi^{\mathrm{Lip}}_P(f)$, it is true that
\begin{equation}
   \ln\left( \E_{\mu}\left[e^{\lambda \cdot f}\right]\right)\, \leq\, \lambda \E_{\mu}[f] + 
    \Theta\left(\frac{\lambda \cdot \Phi^{\mathrm{Lip}}_P(f)}{2 \left(\Lambda^{(1)}_P(\mu)\right)^{1/2}} \right)\, .
\end{equation}
\end{theorem}
This, in turn, gives bounds on $\mathcal{F}^+_{\mu}(f;\cdot)$ using Chebyshev's inequality.
Let us define
\begin{equation}
\label{eq:defXi}
    \Xi_P(\mu,f)\, \stackrel{\mathrm{def}}{:=}\, \frac{2\left( \Lambda^{(1)}_P(\mu)\right)^{1/2}}{\Phi^{\mathrm{Lip}}_P(f)}\, .
\end{equation}
Then we can deduce that the positive fluctuations actually an exponential decay bound almost at the rate of $\Xi_P(\mu,f)$:
\begin{corollary}
\label{cor:AS}
For $a\geq 0$, the positive fluctuations obey the bound (which is written as a negative exponential value for the probability)
\begin{equation}
    -\ln\left(\mathcal{F}^+_{\mu}(f;a)\right)\, \geq\, \max_{0\leq \lambda< \Xi_P(\mu,f)}
    \left(\lambda a - \Theta\left(\frac{\lambda}{\Xi_P(\mu,f)}\right)\right)\, .
\end{equation}
\end{corollary}
By doing a small amount of calculus, this becomes clearer.
\begin{corollary}
\label{cor:CF}
Defining a constant 
\begin{equation}
 \varkappa\, =\, \sum_{n=1}^{\infty} 2^n \ln\left(\frac{1}{1-4^{-n}}\right)\, 
    =\, \sum_{n=1}^{\infty} 2^n \sum_{m=1}^{\infty} \frac{4^{-mn}}{m}\, =\, \sum_{m=1}^{\infty} \frac{2}{m(4^m-2)}\, ,
\end{equation}   
we have 
\begin{equation}
    \Theta(t)\, \leq\, \ln\left(\frac{1}{1-t}\right)+\varkappa\, .    
\end{equation} 
And hence it follows (by using this bound and calculating the Legendre transform of the bounding function) that 
\begin{equation}
    \ln\left(\mathcal{F}^+_{\mu}(f;a)\right)\, \leq\, \varkappa - \vartheta\left(a\cdot \Xi_P(\mu,f)\right)\, ,
\end{equation}
where $\vartheta : [0,\infty) \to [0,\infty)$ is an asymptotically linear function
\begin{equation}
    \vartheta(x)\, =\, x \cdot \varphi(x) + \ln\left(\frac{2}{x}\, \cdot \varphi(x)\right)\, ,\quad 
    \text{ for }\quad \varphi(x)\, =\, \sqrt{1+\frac{1}{x^2}} - \frac{1}{x}\, .
\end{equation}
\end{corollary}
The absolute constant $\varkappa$ satisfies
\begin{equation}
\label{eq:varkappa2}
    1.084640\, \leq\, \varkappa\, \leq\, 1.084645\, .
\end{equation}

An excellent reference for Theorem \ref{thm:AS} and Corollary \ref{cor:AS} is Ledoux's monograph on the concentration of measure phenomenon \cite{Ledoux}.
As stated before, Aida and Stroock proved their result on the way to proving more sophisticated results.
They did not include details of the proofs beyond the basic outline.
But in Section 3.1 in Ledoux's monograph he states the equivalent result as Theorem 3.3, and he gives complete details of the proof.
In particular, when the details are laid out, it becomes apparent that  one can
obtain a slight improvement if one restricts attention to positive fluctuations.
(So far, we have restricted attention to positive fluctuations.)
\begin{theorem}
\label{thm:gen}
In place of (\ref{eq:LipDef}) consider the asymmetric Lipschitz constant
\begin{equation}
\label{eq:aLc}
    \widetilde{\Phi}^{\mathrm{Lip}}_{\mathrm{asym}}(P;f)\, \stackrel{\mathrm{def}}{:=}\, \max_{x \in \mathcal{X}} 
    \left(\sum_{y \in \mathcal{X}} P(x,y) |f(x)-f(y)|^2 \cdot \mathbf{1}_{(0,\infty)}\big(f(x)-f(y)\big)\right)^{1/2}\, .
\end{equation}
Also, define the replacement of (\ref{eq:defXi}) as
\begin{equation}
\label{eq:tildeXiDef}
    \widetilde{\Xi}_P^+(\mu,f)\, \stackrel{\mathrm{def}}{:=}\,  \frac{2\left( \Lambda^{(1)}_P(\mu)\right)^{1/2}}{\widetilde{\Phi}^{\mathrm{Lip}}_{\mathrm{asym}}(P;f)}\, .
\end{equation}
Then, we also have, for 
$\lambda \geq 0$, satisfying $\lambda <  \widetilde{\Xi}_P^+(\mu,f)$, it is true that
\begin{equation}
   \ln\left( \E_{\mu}\left[e^{\lambda \cdot f}\right]\right)\, \leq\, \lambda \E_{\mu}[f] + 
    \Theta\left(\frac{\lambda}{\widetilde{\Xi}_P^+(\mu,f)}\right)\, .
\end{equation}
And, hence, it is also true that
for $a\geq 0$, the positive fluctuations obey the bound 
\begin{equation}
    -\ln\left(\mathcal{F}^+_{\mu}(f;a)\right)\, \geq\, \max_{0\leq \lambda<\widetilde{\Xi}_P^+(\mu,f)}
    \left(\lambda a - \Theta\left(\frac{\lambda}{\widetilde{\Xi}_P^+(\mu,f)}\right)\right)\, .
\end{equation}
\end{theorem}
We do not bother to re-state the calculus facts, but suffice it to say that in Corollary \ref{cor:AS} and Corollary \ref{cor:CF} the constant $\Xi_P(\mu,f)$
may be replaced by $\widetilde{\Xi}_P^+(\mu,f)$.

The slight generalization of Theorem \ref{thm:gen} will be used in our first example. 
We will also include a brief discussion to note that the asymmetric focus on fluctuations is natural and also is already well-established in 
some famous examples.
Let us mention that the asymmetric Lipschitz constant defined in (\ref{eq:aLc}) satisfies good properties.
\begin{proposition} 
\label{prop:PhiProp1}
We have the following properties.
\begin{enumerate}
    \item For all $f,g \in \mathcal{C}(\mathcal{X})$, we have $\widetilde{\Phi}^{\mathrm{Lip}}_{\mathrm{asym}}(P;f+g) \leq \widetilde{\Phi}^{\mathrm{Lip}}_{\mathrm{asym}}(P;f)+\widetilde{\Phi}^{\mathrm{Lip}}_{\mathrm{asym}}(P;g)$.
    \item For all $f \in \mathcal{C}(\mathcal{X})$ and all $c \in [0,\infty)$ we have $\widetilde{\Phi}^{\mathrm{Lip}}_{\mathrm{asym}}(P;c\cdot f) = c \cdot \widetilde{\Phi}^{\mathrm{Lip}}_{\mathrm{asym}}(P;f)$.
    \item Moreover, if we assume that $P$ is irreducible,  then we have:  $\widetilde{\Phi}^{\mathrm{Lip}}_{\mathrm{asym}}(P;f)=0$ implies that $f$ is constant.    
\end{enumerate}
The third condition above is related to being Lipschitz.
Another condition justifying that name is the following.
Suppose that $X_0,X_1,X_2,\dots$ are the states of a random realization of the Markov chain starting from $\mu$, so that
\begin{multline}
    \forall t \in \N\, ,\ \forall x_0,x_1,\dots,x_t \in \mathcal{X}\, , \text{ we have } \\
    \P(\{X_0=x_0\, ,\ X_1=x_1\, , \dots\, ,\ X_t=x_t\})\, =\, \mu(x_0) P(x_0,x_1) \cdots P(x_{t-1},x_t)\, .
\end{multline}
Then
\begin{equation}
   \left( \E\left[\big(f(X_0)-f(X_t)\big)^2 \cdot \mathbf{1}_{(0,\infty)}\big(f(X_0)-f(X_t)\big)\right] \right)^{1/2}\, \leq\, t \cdot \widetilde{\Phi}^{\mathrm{Lip}}_{\mathrm{asym}}(P;f)\, .
\end{equation}
\end{proposition}

\section{First example of the general method: longest increasing sub-sequence of a uniform random permutation}

Given $m,n \in \N = \{1,2,\dots\}$, let $A_m = \left\{0,\frac{1}{m},\frac{2}{m},\dots,1\right\}$ and let $\mathcal{X}_{m,n} = A_m^n$ which is the set of all $(x_1,\dots,x_n)$
with $x_1,\dots,x_n \in A_m$.
The Markov chain transition matrix $P_{m,n} : \mathcal{X}_{m,n} \to \mathcal{X}_{m,n}$
is just replacement of one of the coordinates uniformly at random, where the replacement is by an element of $A_m$ chosen uniformly at random.
Therefore, we have
\begin{equation}
    P_{m,n}\big((x_1,\dots,x_n),(y_1,\dots,y_n)\big)\, =\, \frac{1}{mn}\, \sum_{i=1}^{n} \left(\prod_{j \in \{1,\dots,n\} \setminus \{i\}} \delta(x_i,y_i)\right)\, ,
\end{equation}
where $\delta(x,y) = 1$ if $x=y$, and equals $0$ otherwise.
Since the Kronecker $\delta$ function is symmetric, this is a reversible measure for the invariant measure $\mu_{m,n}$ which is the uniform measure
\begin{equation}
    \mu_{m,n}\big((x_1,\dots,x_n)\big)\, =\, \frac{1}{m^n}\, ,\ \text{ for every }\, (x_1,\dots,x_n) \in A_m^n\, .
\end{equation}
Now suppose that we have a set $J \subset \{1,\dots,n\}$ and we choose for each $i \in J$ a function $f_i : A_m \to \R$ such that $\sum_{x \in A_m} f_i(x) = 0$,
but $f_i$ is not identically zero. In other words, we just assume $f_i$ is orthogonal to the constant function.
Then, letting $F : \mathcal{X}_{m,n} \to \R$ be the function
\begin{equation}
    F(x_1,\dots,x_n)\, =\, \prod_{i \in J} f_i(x_i)\, ,
\end{equation}
it is easy to see that $P_{m,n} F = \frac{n-|J|}{n} F$. (If one of the coordinate indices in $J$ is selected for replacement, then the function $f_i(y_i)$
after replacement has average equal to zero since $f_i$ is orthogonal to the constant function. If any other index is chosen, then the function $F$ is unchanged.)
This suffices to determine a spanning set of eigenvectors.
So the set of eigenvalues is $\{1,1-\frac{1}{n},1-\frac{2}{n},\dots,0\}$.
In particular, using (\ref{eq:DirichletFormula}) and (\ref{eq:GapDef}), we have the following.
\begin{lemma}
For the replacement Markov chain we have been considering, the spectral gap is 
\begin{equation}
    \Lambda_{\mu_{m,n}}^{(1)}(P_{m,n})\, =\, \frac{1}{n}\, .
\end{equation}
\end{lemma}

Now, for the observable, we take $f_{m,n} : \mathcal{X}_{m,n} \to \R$ to be the length of the longest increasing subsequence.
More precisely, let us define $f_n : [0,1]^n \to \R$ to be 
\begin{equation*}
    f_n\big((x_1,\dots,x_n)\big)\, =\, \max\left(\left\{|J|\, :\, J \subset \{1,\dots,n\}\, ,\ \text{ and }\, \forall i,j \in J\, ,\ \text{ we have }\, (i<j) \Rightarrow (x_i<x_j)\right\}\right)\, .
\end{equation*}
Then we take $f_{m,n}$ to just be the restriction $f_n \restriction A_m^n$.

Now let us try to calculate the Lipschitz constant. Suppose that for some $(x_1,\dots,x_n) \in \R^n$, the set $J \subset \{1,\dots,n\}$ is a set where $f_{m,n}\big((x_1,\dots,x_n)\big)=|J|$,
and such that 
\begin{equation}
    \forall i,j \in J\, ,\ \text{ we have }\, (i<j) \Rightarrow (x_i<x_j)\, .
\end{equation}
Then, in one step of the Markov chain, updated by $P_{m,n}$, the only way to decrease the value of $f_{m,n}$
is to choose one of the indices in $J$ to replace by a uniform random sample.
That has probability equal to $|J|/n$.
In that case, it is still possible that the length of the longest increasing subsequence does not decrease.
But if it does decrease, it only decreases by 1.
Therefore, we have, defining
\begin{multline}
    \Delta_{m,n}\big((x_1,\dots,x_n)\big)
    \stackrel{\mathrm{def}}{:=}\,
    \sum_{(y_1,\dots,y_n)\in\mathcal{X}_{m,n}} P_{m,n}\big((x_1,\dots,x_n),(y_1,\dots,y_n)\big) \cdot \\
    \left(f_{m,n}\big((x_1,\dots,x_n)\big)-f_{m,n}\big((y_1,\dots,y_n)\big)\right)^2 \cdot \\
    \mathbf{1}_{(0,\infty)}\left(f_{m,n}\big((x_1,\dots,x_n)\big)-f_{m,n}\big((y_1,\dots,y_n)\big)\right)\,  ,
\end{multline}
it is the case that
\begin{equation}
    \Delta_{m,n}\big((x_1,\dots,x_n)\big)\, \leq\, \frac{|J|}{n}\, .
\end{equation}
Now for the Lipschtiz constant (asymmetric, semi-norm) we have to maximize over all choices of $(x_1,\dots,x_n)$.

That would actually give us a much larger constant than if we restricted to the typical choice of $(x_1,\dots,x_n)$.
That is because the typical value of $|J|$ is approximately $2 \sqrt{n}$, as determined by Vershik and Kerov \cite{VershikKerov}
and Logan and Shepp \cite{LoganShepp}. The correct order is $\sqrt{n}$ from an even easier argument of Hammersley \cite{Hammersley}.
But if we had large deviation bounds, then we could use those to initialize a more refined bound.
In view of all of this, we will just truncate, by-hand.
Given any constant $K \in \N$, let us define
\begin{equation}
    f_n^{(K)}\big((x_1,\dots,x_n)\big)\, =\, \min\left(\left\{K,f_n\big((x_1,\dots,x_n)\big)\right\}\right)\, .
\end{equation}
Then we define $f_{m,n}^{(K)} = f_n^{(K)} \restriction A_m^n$.
Then the above calculations show that if $(x_1,\dots,x_n)$ is such that $|J|\leq K$ then we have
\begin{equation}
    \Delta_{m,n}\big((x_1,\dots,x_n)\big)\, \leq\, \frac{K}{n}\, .
\end{equation}
But if $|J|\geq K+1$ then $f_n\big((x_1,\dots,x_n)\big)=|J|>K$ and no matter what, we will also have $f_n^{(K)}\big((y_1,\dots,y_n)\big)= K= f_n^{(K)}\big((y_1,\dots,y_n)\big)$.
Therefore, from (\ref{eq:aLc})  we have
\begin{equation}
    \widetilde{\Phi}^{\mathrm{Lip}}_{\mathrm{asym}}(P;f^{(K)}_{m,n})\, \leq\, \left(\frac{K}{n}\right)^{1/2}\, . 
\end{equation}
This is the bound which we wanted.

We will take $K$ to be a number depending on $n$, so that we really have a sequence $K_n$.
And we will choose the sequence such that 
\begin{equation}
    \lim_{n \to \infty} n^{-1/2} K_n\, =\, u \in (2,\infty)\, .
\end{equation}
As stated before, to obtain un-restricted bounds, we would need to combine this with large deviation bounds.
Implicitly, we are assuming that in more general applications, it would be easier to get large deviation bounds than concentration-of-measure
bounds.
So, from (\ref{eq:tildeXiDef}) we have
\begin{equation}
    \widetilde{\Xi}_P^+(\mu_{m,n},f_{m,n}^{(K_n)})\, \leq\, \frac{2}{\sqrt{K_n}}\, .
\end{equation}
Then, using this with Corollary \ref{cor:CF}, using $\widetilde{\Xi}_P^+(\mu,f)$ in place of $\Xi_P(\mu,f)$ as discussed at the end of the last section,
we obtain the following.
\begin{corollary}
For the truncation of the length of the longest increasing subsequence,truncated at the level $K=K_n$ such that $\lim_{n \to \infty} K_n/(uN^{1/2}) = 1$
for some $u \in (2,\infty)$, we have the bound
\begin{equation}
    \ln\left(\mathcal{F}^+_{\mu_{m,n}}\left(f_{m,n}^{(K_n)};tn^{1/4})\right)\right)\, \leq\, \varkappa - \vartheta\left(2t \sqrt{\frac{n^{1/2}}{K_n}}\right)\, .
\end{equation}
So in particular, for a fixed $t \in (0,\infty)$ the right hand side converges to $\varkappa - \vartheta(2t/u^{1/2})$. 
\end{corollary}
Note that in the above, the best case for the right hand side would be if one could take $u$ close to $2$ to get $\varkappa - \vartheta(\sqrt{2}\, t)$.
But one cannot do better than that with these methods.
(On can only do that well if one uses good large deviation bounds that are sufficiently good even going down to the median.)
These bounds essentially show that with this method one can determine that the fluctuations are no larger than order $n^{1/4}$.
As shown by Baik, Deift and Johannson the true fluctuations are of order $n^{1/6}$.
But that is bounded by order $n^{1/4}$, so that these bounds are not untrue.
They just are not very sharp.
But that is the situation also for Talagrand's bounds from \cite{Tal}.

We note that the idea of developing asymmetric bounds for the positive and negative fluctuations is not an original idea.
It is already advocated by Talagrand. 
The reason that he obtains bounds as good as he does for the length of the longest increasing subsequence is that the function $f_n\big((x_1,\dots,x_n)\big)$
only depends on $(x_1,\dots,x_n)$ through the points whose indices are in $J \subset \{1,\dots,n\}$.
He called functions such as this ``configuration functions.''
Another good reference is Steele's monograph \cite{Steele}.

We also note that bounds for the negative fluctuations are easily obtained from the bounds for the positive fluctuations, using Markov's inequality.
Of course, one will not obtain as sharp a result that way.
Using 
\begin{equation}
    \int_0^{\infty} \mathcal{F}_{\mu}^{-}(f;t)\, dt\, =\, \int_0^{\infty} \mathcal{F}_{\mu}^{+}(f;t)\, dt\, ,
\end{equation}
we may determine
\begin{equation}
    \mathcal{F}_{\mu}^{-}(f;t)\, \leq\, \frac{1}{t}\, \int_0^{\infty} \mathcal{F}_{\mu}^{+}(f;t)\, dt\, .
\end{equation}
If we have bounds showing that $\mathcal{F}^+_{\mu_{m,n}}\left(f_{m,n}^{(K_n)};tn^{1/4})\right)$ decays exponentially, because we chose $K_n$ to be approximately
$u n^{1/2}$ for some $u \in (2,\infty)$, then we can see that
\begin{equation}
    \int_0^{\infty} \mathcal{F}_{\mu_{m,n}}^{+}(f^{(K_n)}_{m,n};t)\, dt\, \leq Cn^{1/4}\, ,
\end{equation}
for a constant $C$ that depends on $K_n$, or alternatively depends on $u$.
That way, we would see that the negative fluctuations $\mathcal{F}_{\mu_{m,n}}^{-}(f^{(K_n)}_{m,n};t)$ are also decaying when $t$ is at the order of $n^{1/4}$.
So, even though the positive and negative fluctuations have different types of bounds, the 
order of the size of the fluctuations that one obtains bounds for using this technique is the same for both positive and negative fluctuations.
It is order $n^{1/4}$ in this case.

\begin{remark}
If we take $n$ fixed and let $m$ go to $\infty$, then we obtain the analogous bounds when the points $x_1,\dots,x_n$ are chosen uniformly on  the 
continuous interval $[0,1]$,
in an IID fashion.
Since the function only depends on the permutation or relative order induced by the points, that is not a singular limit.
Rather, for finite $m$, the probability that none of the components are equal is 1 minus a quantity which is on the order of $n^2/m$
by the Birthday problem.
When none of the components are equal, conditioning on that event, we do have uniform random permutations,
just as if $x_1,\dots,x_n$ were distributed uniformly on the continuous interval $[0,1]$ in an IID fashion.
\end{remark}

\section{Second example of the general method: Lipschitz functions of 1d Markov chains}

Suppose that $\mathcal{Y}$ is a finite state space, and consider a larger state space $\mathcal{X}_n=\mathcal{Y}^{n+1}$
for some $n \in \N = \{1,\dots,n\}$.
Suppose that $Q:\mathcal{Y}\times \mathcal{Y} \to \R$
is a Markov transition matrix which is irreducible and aperiodic, and suppose that there is a measure $\nu : \mathcal{Y} \to \R$
(satisfying $\nu(x)\geq 0$ for all $x \in \mathcal{Y}$ and $\sum_{x \in \mathcal{Y}}\nu(x)=1$).
And, suppose that $\nu$ is reversible with respect to $Q$:
\begin{equation}
    \forall x,y \in \mathcal{Y}\, ,\ \text{ we have }\, \nu(x) Q(x,y)\, =\, \nu(y) Q(y,x)\, .
\end{equation}
By irreducibility and aperiodicity, we know that $\Lambda^{(1)}_{\nu}(Q)>0$.
The Markov chain we will consider is on $\mathcal{X}_n$ instead of $\mathcal{Y}$.
But this fact is potentially useful for proving lower bounds on the spectral gap of the chain on $\mathcal{X}_n$.

Before stating the Markov transition matrix for the chain on $\mathcal{X}_n$, let us define the measure we wish to be the invariant measure for 
the Markov chain.
Let $\mu_n : \mathcal{X}_n \to \R$ be the measure defined as 
\begin{equation}
    \mu_n\big((x_0,x_1,\dots,x_n)\big)\, =\, \mathbb{P}_{\nu}(\{X_0=x_0\, ,\ X_1=x_1\, ,\ \dots\, ,\ X_n=x_n\})\, =\, \nu(x_0) P(x_0,x_1) \cdots P(x_{n-1},x_n)\, ,
\end{equation}
for each $(x_0,x_1,\dots,x_n) \in \mathcal{X}_n = \mathcal{Y}^{n+1}$.
Here $X_0,X_1,\dots,X_n$ is viewed as a sample of the 1d Markov chain, from times $t=0$ to $t=n$, started at time $0$ in the distribution $\nu$.
By reversiblity of $\nu$ with respect to $Q$, we also have
\begin{equation}
    \mu_n\big((x_0,x_1,\dots,x_n)\big)\, =\, \nu(x_n) P(x_n,x_{n-1}) \cdots P(x_{1},x_0)\, ,
\end{equation}
and for $t \in \{1,\dots,n-1\}$
\begin{equation}
    \mu_n\big((x_0,x_1,\dots,x_n)\big)\, =\, \nu(x_t) \Big(P(x_t,x_{t+1}) \cdots P(x_{n-1},x_n)\Big) \Big(P(x_t,x_{t-1})\cdots P(x_1,x_0)\Big)\, .
\end{equation}
These alternative formulations are potentially useful for proving reversibility for the Markov chain on $\mathcal{X}_n$, which is what we consider next.

We consider Glauber dynamics as the Markov chain on $\mathcal{X}_n=\mathcal{Y}^{n+1}$.
In other words, we consider $P_n : \mathcal{X}_n \times \mathcal{X} \to \R$ to be the Markov transition matrix where
\begin{equation}
    P_n\big((x_0,\dots,x_n),(y_0,\dots,y_n)\big)\, =\, \frac{1}{n+1}\, \sum_{t=0}^n P_{n,t}\big((x_0,\dots,x_n),(y_0,\dots,y_n)\big)\, ,
\end{equation}
where for $t \in \{1,\dots,n-1\}$ we have
\begin{equation}
    P_{n,t}\big((x_0,\dots,x_n),(y_0,\dots,y_n)\big)\, =\, \frac{Q(x_{t-1},y_t)Q(y_t,x_{t+1})}{\sum_{z\in\mathcal{Y}}Q(x_{t-1},z)Q(z,x_{t+1})}\, \cdot
    \prod_{s \in \{0,\dots,n\} \setminus \{t\}} \delta(x_s,y_s)\, ,
\end{equation}
while
\begin{equation}
    P_{n,0}\big((x_0,\dots,x_n),(y_0,\dots,y_n)\big)\, =\, \frac{\nu(y_0) Q(y_0,x_1)}{\sum_{z\in\mathcal{Y}}\mu(z) Q(z,x_1)}\, \cdot
    \prod_{s \in \{1,\dots,n\}} \delta(x_s,y_s)\, .
\end{equation}
Note that by reversibility 0f $\nu$ with respect to $Q$, these can be written seemingly different but equivalent ways.
For example,
\begin{equation}
    P_{n,0}\big((x_0,\dots,x_n),(y_0,\dots,y_n)\big)\cdot \mu\big((x_0,x_1,\dots,x_n)\big)\,\, =\, \frac{Q(x_1,y_0)}{\sum_{z\in\mathcal{Y}} Q(x_1,z)}\, \cdot
    \prod_{s \in \{1,\dots,n\}} \delta(x_s,y_s)\, ,
\end{equation}
as well.

It is easy to see that each of the $P_{n,t}$ matrices, for $t \in \{0,1,\dots,n-1\}$ is such that $\mu_n$ is reversible for $P_{n,t}$.
For example, if $t$ is in $\{1,\dots,n-1\}$, then
\begin{equation}
    P_{n,t}\big((x_0,\dots,x_n),(y_0,\dots,y_n)\big)\, =\, \frac{Q(x_{t-1},y_t)Q(y_t,x_{t+1})}{\sum_{z\in\mathcal{Y}}Q(x_{t-1},z)Q(z,x_{t+1})}\, \cdot 
    \nu(x_0) Q(x_0,x_1)\cdots Q(x_{n-1},x_n)\, ,
\end{equation}
if $y_s=x_s$ for all $s \in \{1,\dots,n\} \setminus \{t\}$ (and it equals $0$ otherwise).
Isolating $x_t$ and $y_t$, and assuming $y_s=x_s$ for all $s \in \{1,\dots,n\} \setminus \{t\}$ in order not to get 0, 
this is 
\begin{multline}
\label{eq:RevSym}
    P_{n,t}\big((x_0,\dots,x_n),(y_0,\dots,y_n)\big) \mu_n\big((x_0,x_1,\dots,x_n)\big)\, =\, Q(x_{t-1},y_t)Q(y_t,x_{t+1})Q(x_{t-1},x_t)Q(x_t,x_{t+1}) \\
    \cdot F_{n,t}(x_0,\dots,x_{t-1},x_{t+1},\dots,x_n)\, ,
\end{multline}
for
\begin{equation}
    F_{n,t}(x_0,\dots,x_{t-1},x_{t+1},\dots,x_n)\, =\, \frac{\nu(x_0) Q(x_0,x_1) \cdots Q(x_{t-2},x_{t-1}) \cdot Q(x_{t+1},x_{t+2}) \cdots Q(x_{n-1},x_n)}{\sum_{z\in\mathcal{Y}}Q(x_{t-1},z)Q(z,x_{t+1})}\, .
\end{equation}
Clearly (\ref{eq:RevSym}) is symmetric in interchange of the two coordinates of $(x_t,y_t)$.
Also, the conditions imposed by multiplying by $\prod_{s \in \{0,\dots,n\} \setminus \{t\}} \delta(x_s,y_s)$ is also symmetric in interchange of every $(x_s,y_s)$
for $s \in \{0,\dots,n\} \setminus \{t\}$.
So $\mu$ is a reversible measure for Glauber dynamics.
We refer to Chapter 3 of Levin and Peres, for example \cite{LevinPeres}, for more details on Glauber dynamics.

\begin{proposition}
For the Glauber dynamics we have been considering, there is a constant $\lambda_1$ satisfying
\begin{equation}
    \lambda_1\, >\, 0\, ,
\end{equation}
such that
\begin{equation}
    \Lambda^{(1)}_{\mu_n}(P_n)\, \geq\, \frac{\lambda_1}{n+1}\, .
\end{equation}
\end{proposition}
We will not prove this, here.
But it is reportedly well-known. A reference is Lu and Yau's paper on Glauber dynamics and Kawasaki dynamics \cite{LuYau}.

We have a specific application in mind, which we call target frequency analysis, which is also hypothesis testing for the power spectrum (Fourier transform amplitudes-squared) integrated over certain intervals.
But before moving to that example, let us just quickly state the general result.
\begin{corollary}
\label{cor:Glauber}
Suppose that we have a function $f_n : \mathcal{X}_n \to \R$ satisfying 
\begin{equation}
    \widetilde{\Phi}^{\mathrm{Lip}}_{\mathrm{asym}}(P_n;f_n)\, \leq\, \frac{C}{(n+1)^p}\, ,
\end{equation}
for some power $p$. Then we have the bound
\begin{equation}
    \ln\left(\mathcal{F}^+_{\mu_n}\left(f_n;a (n+1)^{-p+\frac{1}{2}}\right)\right)\, \leq\, \varkappa - \vartheta\left(a\cdot \frac{2\sqrt{\lambda_1}}{C}\right)\, ,
\end{equation}
using the notation from Corollary \ref{cor:CF}.
\end{corollary}

For us, the power we will obtain will be $p=1$, so that the fluctuations will be shown to be bounded by $O(n^{-1/2})$ in this way, for a nonnegative observable
whose mean is order-1.

\subsection{Example: Target frequency analysis for Markov chains}
As another basic application, we consider a statistic for time series which was considered by the authors and Jung in \cite{FJS}.
This is called ``target frequency analysis.''
The application is important in biostatistics. 
But it also supplies a pedagogically valuable example for the technique.

Let $\mathcal{Y}$ be $\mathcal{Y}_{m,\epsilon} = \{-m\epsilon,(-m+1)\epsilon,\dots,m\epsilon\}$.
For us, an important quantity is the radius of this chain $R = m\epsilon$.
Note that $\mathcal{Y}_{m,\epsilon} \subset \R$, and the next step in the description of the function $f_{n}$ on $\mathcal{X}_n=(\mathcal{Y}_{m,\epsilon})^{n+1}$
only relies on that.
Given a real sequence $(x_0,\dots,x_n)$ we define the Fourier transform
$\phi_{(x_0,\dots,x_n)} : \Z \to \C$ defined as if $(x_0,\dots,x_n)$ were the $n+1$ components of a periodic signal
\begin{equation*}
    \phi_{(x_0,\dots,x_n)}(k)\, =\, \frac{1}{\sqrt{n+1}}\, \sum_{t=0}^n \exp\left(\frac{-2\pi i k t}{n+1}\right) x_t\, ,
\end{equation*}
where $i = \sqrt{-1}$, as usual (despite the fact that in earlier sections the symbol $i$ was used for an integer index).
The choice of the prefactor $1/\sqrt{n+1}$ is such that Parseval's identity is satisfied
\begin{equation}
\label{eq:Pars}
\begin{split}
    \sum_{k=0}^{n} |\phi_{(x_0,\dots,x_n)}(k)|^2\, 
    &=\, \frac{1}{n+1} \sum_{k=0}^{n} \sum_{s=0}^{n} \sum_{t=0}^n \exp\left(\frac{-2\pi i k (t-s)}{n+1}\right) x_tx_s\\
    &=\, \frac{1}{n+1} \sum_{s=0}^{n} \sum_{t=0}^n \left(\sum_{k=0}^{n} \exp\left(\frac{-2\pi i k (t-s)}{n+1}\right)\right) x_tx_s\\
    &=\, \frac{1}{n+1} \sum_{s=0}^{n} \sum_{t=0}^n \big((n+1) \delta(s,t)\big) x_tx_s\\
    &=\, \sum_{t=0}^n x_t^2\, .
\end{split}
\end{equation}
Now, given any choice of $a,b \in \Z$ satisfying $0<a<b<(n+1)/2$, we consider the observable of interest to be
\begin{equation}
    f_n(x_0,\dots,x_n)\, =\, \frac{1}{n+1} \sum_{k=a}^b |\phi_{(x_0,\dots,x_n)}(k)|^2\, .
\end{equation}
In other words, using the language of signal processing, it is the power spectrum integrated from $a$ to $b$.
Now it is rescaled by $1/(n+1)$ because for a signal of length $n+1$, we expect the total $\ell^2$-norm (also called the total power)
to be of order $(n+1)$.
So this rescales to give an order-1 quantity.
Note that the Fourier transform is an isometry by (\ref{eq:Pars}), therefore $f_n$ may be viewed as a contraction mapping times a constant $1/(n+1)$.
For this reason, we obtain
\begin{equation}
    \Phi^{\mathrm{Lip}}_{P_n}(f_n)\, \leq\, \frac{2R^2}{n+1}\, .
\end{equation}
Since 
$\widetilde{\Phi}^{\mathrm{Lip}}_{\mathrm{asym}}(P_n;f_n)\leq \Phi^{\mathrm{Lip}}_{P_n}(f_n)$, this proves the following
using Corollary \ref{cor:Glauber}.
\begin{corollary}
For the function $f_n$ written above, we have
\begin{equation}
    \ln\left(\mathcal{F}^+_{\mu_n}\left(f_n;\frac{a}{\sqrt{n+1}}\right)\right)\, \leq\, \varkappa - \vartheta\left(a\cdot \frac{\sqrt{\lambda_1}}{R^2}\right)\, ,
\end{equation}
\end{corollary}
We call the function $f_n$ by the name ``target frequency analysis.''
It has a special property: if we replace the present set-up by the case where $\mathcal{Y}=\R$ and allow $X_0,\dots,X_n$
to be IID standard, normal random variables (also called white noise, by some),
then the Fourier transform has the property for frequencies $k$ satisfying $0<k<(n+1)/2$ the real and imaginary parts are all IID standard, normal
random variables. (Here IID refers to the independence of the real and imaginary parts, as well as independence for different values of $k$.)
By the Parseval identity, isometry property, it is elementary that the Fourier transform of IID complex-valued signals with real and imaginary
parts being IID standard, normal random variables would have the same property.
But the property stated for the Fourier transform of a real signal is slightly less trivial, although it may be easily checked using covariance matrices.
One may also thinking of this fact as arising from the slight extra information included in the dihedral symmetry over the usual cyclic symmetry, for the dihedral group
$D_n$ being the semi-direct product of the cyclic group $(\mathbb{Z}/n\mathbb{Z})$ with the involution group $(\mathbb{Z}/2\mathbb{Z})$.
One can also see it by using properties of the complex conjugation, which amounts to the same.
But it is probably the simplest example of a more general phenomenon where for special symmetrical models, random variables defined 
on a large space have unexpected projections into some components which also have simple, explicit distributions on smaller spaces.

The corollary above shows that more generally, for Markov chain models of a time series, the {\em target frequency analysis}
will still be concentrating at least in the sense that the fluctuations are no larger than order $n^{1/2}$.


\baselineskip=12pt
\bibliographystyle{plain}

\noindent
\underline{Contact:}
\texttt{slstarr@uab.edu}

\end{document}